\documentclass[12]{amsart}
  
\usepackage{layout}
\usepackage{amsfonts}
\usepackage{amssymb}
\usepackage{amsthm}
\newtheorem{theorem}{Theorem}
\newtheorem{thmx}{Theorem}

\def\r{\mathbb R}
\def\l{\mathbb L}

\begin{document}

\title{Translation surfaces of linear Weingarten type}

\author{Antonio Bueno and Rafael L\'opez}
\address{Departamento de Geometr\'{\i}a y Topolog\'{\i}a\\
Universidad de Granada\\
18071 Granada. Spain}
\email{rcamino@ugr.es}
\email{buenolinares@gmail.com}
\thanks{Partially supported by MEC-FEDER
 grant no. MTM2011-22547 and Junta de Andaluc\'{\i}a grant no. P09-FQM-5088.}

\begin{abstract}We give a relatively simple proof  that a translation surface in Euclidean space that satisfies a relation of type $aH+bK=c$, for some real numbers $a,b,c$, where $H$ and $K$ are the mean curvature and the Gauss curvature of the surface, respectively, must have $a=0$ or $b=0$, and thus,  $K$ is constant or $H$ is constant. Our method of proof extends to the Lorentzian ambient space.
\end{abstract}


\maketitle
\section{Introduction and results.}

A Weingarten surface in Euclidean space $\r^3$ is a surface $S$ whose mean curvature $H$ and Gauss curvature $K$ satisfies a non-trivial relation $\Psi(H,K)=0$. This type of surfaces were introduced by the very Weingarten in the context of the problem of finding all surfaces isometric to a given surface of revolution and have been extensively studied in the literature \cite{we}. In order to simplify the study of Weingarten surfaces, it is natural to impose some added geometric condition on the surface, as for example, that $S$ is ruled or rotational \cite{be,dk,di,ku,rs}.

Following this strategy,  Dillen, Goemans and Van de Woestyne considered Weingarten surfaces that are graphs of type  $z=f(x)+g(y)$, where $f$ and $g$ are smooth functions defined in some intervals $I,J\subset\r$, respectively \cite{dgv}. A surface $S$ in $\r^3$ is called a {\it translation surface} if it can locally parametrize as $X(x,y)=(x,y,f(x)+g(y))$. In particular,  a translation surface $S$ has the property that the translations of a parametric curve $x=ct$ by the parametric curves $y=ct$ remain in $S$ (similarly for the parametric curves $x=ct$). In the cited paper, the authors classify all translation surfaces of Weingarten type:

\begin{thmx}[\cite{dgv}]\label{tdgv}
A translation surface in $\r^3$ of Weingarten type is a plane, a generalized cylinder, a Scherk's minimal surface or an elliptic paraboloid.
\end{thmx}

The proof given in \cite{dgv} (see also \cite{go})  discusses many cases and it involves the solvability of a large number of ODE systems. In fact,   in \cite{dgv} it is described the procedure and it requires of calculations which are done with a computer  program (as Maple) to manipulate the algebraic operations. This is the reason that  some authors previously obtained  partial results assuming simpler functions $f$ and $g$, as for example, that they are polynomial in its variables, simplifying and doing easier the computations (\cite{mn,wo}).

In this paper  we provide  a significantly simpler proof of Th. \ref{tdgv} when the Weingarten relation is linear in its variables. A \emph{linear Weingarten surface}  in Euclidean space $\r^3$ is a surface where there exists a relation
\begin{equation}\label{w1}
a\ H+b\ K=c,
\end{equation}
for some real numbers $a,b,c$, not all zero.   In the class of linear Weingarten surfaces, we mention  two families of surfaces that correspond
with trivial choices of the constants $a$ and $b$: surfaces with
constant Gauss curvature ($a=0$) and surfaces
with constant mean curvature ($b=0$).   In Th. \ref{tdgv},  only the three first surfaces are linear Weingarten surfaces, which have constant $H$ or constant $K$:  a plane ($H=K=0$), a generalized cylinder ($K=0$) and  the Scherk's minimal surface parametrized as
$z=\log(\cos (\lambda y))-\log(\cos(\lambda x))$, $\lambda>0$ ($H=0$).  Besides these two families of surfaces, the  classification of linear Weingarten surfaces in the general case is almost completely open today. See \cite{gmm,lo1,rs}.

The result that we prove is:

\begin{theorem}\label{t1}
A translation surface in Euclidean space $\r^3$ of linear Weingarten type is a surface with constant Gauss curvature $K$ or constant mean curvature  $H$. In particular, the surface is congruent with a  plane,  a generalized cylinder or a Scherk's minimal surface.
\end{theorem}

This proves that in the family of translation surfaces, there doesn't exist new  linear Weingarten surfaces besides the trivial choices of $a,b$ in \eqref{w1}. We point out that an early work of Liu proved    that the only translations surfaces with constant $K$ or constant $H$  are the three first surfaces of  Th. \ref{t1} (\cite{li}). Finally, and with minor modifications, we extend in Th. \ref{t2} our results to the Lorentzian ambient space (see also \cite{dgv}).

\section{Proof of Theorem \ref{t1}}
 The mean curvature $H$ and the Gauss curvature $K$ are expressed in a local parametrization $X$ as
\begin{equation}\label{hk1}
H=\frac{eG-2fF+gE}{2(EG-F^2)},\ \ K=\frac{eg-f^2}{EG-F^2},
\end{equation}
where $\{E,F,G\}$ and $\{e,f,g\}$ are the coefficients of the first fundamental form and the second fundamental form, respectively. Assume that $S$ is a translation surface expressed locally as $X(x,y)=(x,y,f(x)+g(y))$ for some smooth functions $f$ and $g$.  Then $H$ and $K$ are
\begin{equation}\label{hk2}
H=\frac{f''(1+g'^2)+g''(1+f'^2)}{2(1+f'^2+g'^2)^{\frac{3}{2}}},\ \ K=\frac{f''g''}{(1+f'^2+g'^2)^2}.
\end{equation}
Suppose now that $S$ is also a linear Weingarten surface, where $H$ and $K$ satisfy the linear relation \eqref{w1}. The proof of Theorem \ref{t1} is by contradiction and we suppose that $a,b\not=0$. Let us observe that this implies  $f''\not=0$ and  $g''\not=0$ because on the contrary, and from \eqref{hk2},  $H$ is constant. Let
$$W=EG-F^2=1+f'^2+g'^2.$$
We distinguish two cases according the value of $c$.

\subsection{Case $c=0$.}

Suppose $c=0$ in \eqref{w1}.  With the change $a\rightarrow 2a$ and by using \eqref{hk2}, Equation \eqref{w1} writes as
\begin{equation}\label{c0}
a\ \frac{f''(1+g'^2)+g''(1+f'^2)}{(1+f'^2+g'^2)^{\frac{3}{2}}}+b\ \frac{f''g''}{(1+f'^2+g'^2)^2}=0.
\end{equation}
We  multiply \eqref{c0} by $W^2$ and divide by $(1+g'^2)(1+f'^2)$ obtaining
\begin{equation}\label{c1}
a\ \left(\frac{f''}{1+f'^2}+\frac{g''}{1+g'^2}\right)\sqrt{W}+b\ \frac{f''}{1+f'^2}\frac{g''}{1+g'^2}=0.
\end{equation}
Introduce the next notation:
\begin{equation}\label{efg}
F=\frac{f''}{1+f'^2},\ \ G=\frac{g''}{1+g'^2}.
\end{equation}
In particular, since $f''\not=0$ and $g''\not=0$, then $F\not=0$ and $G\not=0$. Then \eqref{c1} writes as
\begin{equation}\label{c2}
a(F+G)\sqrt{W}+bFG=0.
\end{equation}
Let us observe that this identity implies $F+G\not=0$, since on the contrary, $b=0$. From \eqref{c2}, we have
$$1+f'^2+g'^2=W=\frac{b^2}{a^2}\left(\frac{FG}{F+G}\right)^2.$$
We differentiate this equation with respect to $x$ and next, with respect to $y$. Because the left hand side is a sum of a function of $x$ and a function $y$, this calculation yields $0$. On the other hand,  the right hand side concludes
\begin{equation}\label{eab}
6\frac{b^2}{a^2}\frac{F^2G^2F'G'}{(F+G)^4}=0.
\end{equation}
This implies $F'=G'=0$ and thus, $F$ and $G$ are constants. From \eqref{c2}, we deduce that $W=1+f'^2+g'^2$ is constant, in particular, $f'$ and $g'$ are constant: a contradiction with the fact that $f'',g''\not=0$.

\subsection{Case $c\not=0$.}

Consider $c\not=0$ in \eqref{w1}. Dividing by $c$, and after a change of notation, the relation \eqref{w1} writes as
\begin{equation}\label{nc0}
a\ \frac{f''(1+g'^2)+g''(1+f'^2)}{(1+f'^2+g'^2)^{\frac{3}{2}}}+b\ \frac{f''g''}{(1+f'^2+g'^2)^2}=1,
\end{equation}
or equivalently
\begin{equation}\label{nc2}
a(F+G)\sqrt{W}+bFG=\frac{W^2}{(1+f'^2)(1+g'^2)},
\end{equation}
where $F$ and $G$ are given in \eqref{efg}.
We differentiate \eqref{nc2} separately with respect to $x$ and with respect to $y$:
$$a\left(F'\sqrt{W}+(F+G)\frac{f'f''}{\sqrt{W}}\right)+b F'G=\frac{4Wf'f''}{(1+f'^2)(1+g'^2)}-\frac{2f'f''W^2}{(1+f'^2)^2(1+g'^2)}.$$
$$a\left(G'\sqrt{W}+(F+G)\frac{g'g''}{\sqrt{W}}\right)+b F'G=\frac{4Wg'g''}{(1+f'^2)(1+g'^2)}-\frac{2g'g''W^2}{(1+f'^2)(1+g'^2)^2}.$$
Dividing the first equation by $f'f''$ and the second one by $g'g''$, we have
$$a\frac{F'\sqrt{W}}{f'f''}+b\frac{F'G}{f'f''}+\frac{2W^2}{(1+f'^2)^2(1+g'^2)}=
a\frac{G'\sqrt{W}}{g'g''}+b\frac{FG'}{g'g''}+\frac{2W^2}{(1+f'^2)(1+g'^2)^2}.$$
From \eqref{nc2}, we replace the value of $W^2$ in the above expression, obtaining
\begin{eqnarray}\label{nc3}
&&a\left(\frac{F'}{f'f''}+\frac{2(F+G)}{1+f'^2}-\frac{G'}{g'g''}-\frac{2(F+G)}{1+g'^2}\right)\sqrt{W}\nonumber \\
&&+
b\left(\frac{F'G}{f'f''}+\frac{2FG}{1+f'^2}-\frac{FG'}{g'g''}-\frac{2FG}{1+g'^2}\right)=0.
\end{eqnarray}
Now we write \eqref{nc0} as
$$a\left(f''(1+g'^2)+g''(1+f'^2)\right)\sqrt{W}+bf''g''=W^2$$
and we differentiate this expression with respect to $x$ and with respect to $y$:
\begin{eqnarray*}
&&a\left(f'''(1+g'^2)+2f'f''g''\right)\sqrt{W}+a\left(f''(1+g'^2)+g''(1+f'^2)\right)\frac{f'f''}{\sqrt{W}}\\
&& +bf'''g''=4f'f''W.\end{eqnarray*}
\begin{eqnarray*}
&& a\left(2f''g'g''+g'''(1+f'^2)\right)\sqrt{W}+a\left(f''(1+g'^2)+g''(1+f'^2)\right)\frac{g'g''}{\sqrt{W}}\\
&&+bf''g'''=4g'g''W.\end{eqnarray*}
From both equations, we obtain the value of $W$ on the right hand sides and we equal both expressions, deducing
\begin{equation}\label{nc5}
a\left(\frac{f'''}{f'f''}(1+g'^2)+2g''-2f''-\frac{g'''}{g'g''}(1+f'^2)\right)\sqrt{W}=b\left(f''\frac{g'''}{g'g''}-g''\frac{f'''}{f'f''}\right).
\end{equation}
If we write \eqref{nc3} and \eqref{nc5} as $P_1\sqrt{W}=Q_1$ and $P_2\sqrt{W}=Q_2$, respectively, we obtain $P_1Q_2-P_2Q_1=0$. After some manipulations, this identity writes as
$$\left(f'f''^2g'''-f'''g'g''^2\right)\left(2f'f''g'g''(f''-g'')+f'f''(1+f'^2)g'''-f'''g'g''(1+g'^2)\right)=0,$$
that is, $P_2Q_2=0$. We discuss by cases:
\begin{enumerate}
\item Case $P_2=0$ and $Q_2\not=0$. Then \eqref{nc5} implies $a=0$, a contradiction.
\item Case $P_2\not=0$ and $Q_2=0$ Then \eqref{nc5} implies $b=0$, a contradiction.
\item Case $P_2=Q_2=0$. These two equations write as
\begin{equation}\label{nc61}\frac{f'''}{f'f''^2}=\frac{g'''}{g'g''^2}\end{equation}
\begin{equation}\label{nc62}
2(f''-g'')+\frac{g'''}{g'g''}(1+f'^2)-\frac{f'''}{f'f''}(1+g'^2)=0.
\end{equation}
Equation \eqref{nc61} implies the existence of $\lambda\in\r$ such that
\begin{equation}\label{nc7}
\frac{f'''}{f'f''^2}=\frac{g'''}{g'g''^2}=2\lambda
\end{equation}
and thus
$$\frac{f'''}{f'f''}=2\lambda f'',\ \ \frac{g'''}{g'g''}=2\lambda g''.$$
Substituting the above in \eqref{nc62}, we get
$$
2(f''-g'')+2\lambda(1+f'^2) g''-2\lambda (1+g'^2)f''=0,
$$
or
\begin{equation}\label{nc8}
f''-g''+\lambda g''-\lambda f''=\lambda f''g'^2-\lambda g''f'^2.
\end{equation}
If $\lambda\neq 0$, differentiating this equation with respect to $x$ and then with respect to  $y$, we deduce
$$f'f''g'''=g'g''f'''.$$
As we suppose that $f'',g''\neq 0$, we conclude that
$$\frac{f'''}{f'f''}=\frac{g'''}{g'g''}=\mu$$
for some constant $\mu\in\r$.  Substituting in \eqref{nc7} we deduce that $f'',g''$ are both constant functions, so \eqref{nc7} yields to $\lambda$ being zero, a contradiction.

Therefore,  $\lambda=0$ in \eqref{nc7}. Equation \eqref{nc8} says now that $f''=g''=m$, for some real number $m\neq 0$. Then \eqref{nc0} writes as
$$am(2+f'^2+g'^2)=W^{\frac32}-bm^2W^{-\frac12}.$$
Differentiating with respect to $x$ and simplifying by $f'f''$,  we get
$$2am=3W^{\frac{1}{2}}+bm^2 W^{-\frac{3}{2}},$$
which implies that $W$ is constant and this would say that $f''=g''=0$, a contradiction.
\end{enumerate}

\section{The Lorentzian case}\label{s-lore}

We consider the Lorentz-Minkowski space $\l^3$, that is, the real vector space $\r^3$ endowed with the metric $(dx)^2+(dy)^2-(dz)^2$ where $(x,y,z)$ are the canonical coordinates. A surface $S$ immersed in $\l^3$ is said non degenerate if the induced metric on $S$ is non degenerated. The induced metric on $S$ can only be of two types: positive definite and the surface is called spacelike, or a Lorentzian metric, and the surface is called timelike. For both types of surfaces, it is defined the mean curvature $H$ and the Gauss curvature $K$ and we say again that the surface is of linear Weingarten type if there exists a linear relation between $H$ and $K$ as in \eqref{w1}.

Similarly, in Lorentzian setting we can extend the concept of translation surface. A  surface $S$ in $\l^3$ is again locally a graph on one of the coordinate planes, since this property is not metric but because $S$ is immersed. Thus a translation surface in $\l^3$ is a surface that writes locally as the graph of a function which is the sum of two real functions. However, in $\l^3$  we can say a bit more.  If $S$ is spacelike, then $S$ is a graph on the $xy$-plane and if $S$ is a timelike surface, then $S$ is a graph on the $xz$-plane or on the $yz$-plane \cite{wei}. Therefore, if $S$ is a translation surface in $\l^3$, we may suppose that:
\begin{enumerate}
\item If $S$ is spacelike, then $S$ writes  locally as $z=f(x)+g(y)$.
 \item If $S$ is timelike, then $S$ writes locally as  $y=f(x)+g(z)$ or as $x=f(y)+g(z)$.
 \end{enumerate}
In \cite{dgv}, Theorem \ref{tdgv} was extended to non-degenerate surfaces of $\l^3$, obtaining a similar result. Again, in this classification, the only translation surfaces of linear Weingarten type appear with trivial choices of $a$ and $b$ and  the surfaces have constant $H$ or constant $K$. Similarly, we extend Theorem \ref{t1} as follows:

\begin{theorem}\label{t2}
A translation non-degenerate surface in Lorentz-Minkowski space $\l^3$ of linear Weingarten type is a surface with constant Gauss curvature $K$ or constant mean curvature  $H$.
\end{theorem}

Translations surfaces in $\l^3$ with constant mean curvature or constant Gauss curvature were classified in \cite{li} and they are a plane, a Scherk's minimal surface or a generalized cylinder.

\begin{proof} The proof of Th. \ref{t2} is similar as Th. \ref{t1} and we only sketch the differences. Moreover, we will carry jointly the cases that the surface $S$ is spacelike or timelike. Again, we suppose by contradiction that $a,b\not=0$ in \eqref{w1}. The  expressions of $H$ and $K$ in local coordinates are
$$H=\epsilon\ \frac12\frac{eG-2fF+gE}{EG-F^2},\ \ K= \epsilon\ \frac{eg-f^2}{EG-F^2},$$
where $\epsilon=-1$ is $S$ is spacelike and $\epsilon=1$ if $S$ is timelike (\cite{lo2,wei}).

Suppose that $S$ writes as $z=f(x)+g(y)$ if $S$ is spacelike or  $y=f(x)+g(z)$ if $S$ is timelike. Then
$$H=\epsilon\frac{-\epsilon f''(1-g'^2)+g''(1+\epsilon f'^2)}{2( (1+\epsilon f'^2-g'^2))^{\frac{3}{2}}}, \ K=-\frac{f''g''}{(1+\epsilon f'^2-g'^2)^2},$$
with  $W=1+\epsilon f'^2-g'^2>0$. Let
$$F=\frac{f''}{1+\epsilon f'^2},\ \ G=\epsilon \frac{g''}{-1+g'^2}.$$
If $c=0$ in \eqref{w1}, then \eqref{c2} is the same, obtaining \eqref{eab}. This implies that $W$ is constant, a contradiction.

If $c\not=0$, then we assume after a change of constants $a$ and $b$ that $c=1$. Now the linear Weingarten  condition \eqref{w1} expresses as
\begin{equation}\label{lnc2}
a(F+G)\sqrt{W}+bFG=\varepsilon\frac{W^2}{(1+\epsilon f'^2)(-1+g'^2)}.
\end{equation}
Now \eqref{nc3} and \eqref{nc5} write, respectively, as
\begin{eqnarray*}
& &a\left(\frac{F'}{f'f''}+\frac{2(F+G)}{\varepsilon+f'^2}+\varepsilon\frac{G'}{g'g''}+\varepsilon\frac{2(F+G)}{-1+g'^2}\right)\sqrt{W}\\
& & +
b\left(\frac{F'G}{f'f''}+\frac{2FG}{\varepsilon+f'^2}+\varepsilon\frac{FG'}{g'g''}+\varepsilon\frac{2FG}{-1+g'^2}\right)=0\end{eqnarray*}
\begin{eqnarray*}
& &a\left(\frac{f'''}{f'f''}(-1+g'^2)+2g''+2\varepsilon f''+\varepsilon(\varepsilon+f'^2)\frac{g'''}{g'g''}\right)\sqrt{ W}\\
& & +\varepsilon b\left(f''\frac{g'''}{g'g''}+\varepsilon g''\frac{f'''}{f'f''}\right)=0.
\end{eqnarray*}
We deduce
\begin{eqnarray*}
&& \left(f'f''^2g'''+\varepsilon f'''g'g''^2\right)\left(2f'f''g'g''(f''+\varepsilon g'')+f'f''(f'^2+\varepsilon)g'''\right.\\
&& \left.+\varepsilon f'''g'g''(g'^2-1)\right)=0
\end{eqnarray*}
and now the discussion by cases is similar  as it was done in the Euclidean case, obtaining that $W$ is constant, a contradiction.

\end{proof}


\end{document}